\documentclass[12pt]{article}
\usepackage{a4wide}
\usepackage{amsfonts}
\usepackage{amsmath}
\usepackage{amssymb}
\usepackage{mathrsfs}
\usepackage[all]{xy}
\newtheorem{definition}{Definition}[section]
\newtheorem{lemma}[definition]{Lemma}
\newtheorem{theorem}[definition]{Theorem}
\newtheorem{proposition}[definition]{Proposition}
\newtheorem{corollary}[definition]{Corollary}
\newtheorem{Remark}[definition]{Remark}
\newenvironment{remark}{\begin{Remark}\rm}{\end{Remark}}
\newcommand{\cst}{\mathrm{C}^*\!}

\newcommand{\CC}{{\mathbb C}}
\newcommand{\NN}{{\mathbb N}}
\newcommand{\ZZ}{{\mathbb Z}}
\newcommand{\TT}{{\mathbb T}}
\newcommand{\FF}{{\mathbb F}}
\newcommand{\id}{\mathrm{id}}
\newcommand{\comp}{\!\circ\!}
\newcommand{\bez}{\setminus}

\newcommand{\algt}{\otimes_{\scriptscriptstyle\mathrm{alg}}\!}
\newcommand{\its}[3]{\left(#1\,\vline\,#2\,\vline\,#3\right)}
\newcommand{\K}[1]{{\mathcal K}\!\left(#1\right)}
\newcommand{\B}[1]{B\!\left(#1\right)}
\newcommand{\eps}{\varepsilon}
\newcommand{\ph}{\varphi}
\newcommand{\del}{\delta}
\newcommand{\tens}{\otimes}
\newcommand{\Mor}[2]{\mathrm{Mor}\!\left(#1,#2\right)}
\newcommand{\M}[1]{\mathrm{M}\!\left(#1\right)}
\newcommand{\Rep}[2]{\mathrm{Rep}\!\left(#1,#2\right)}
\newcommand{\Aff}[1]{#1^{\eta}}
\newcommand{\Ci}[1]{\mathrm{C}_{\infty}\!\left(#1\right)}
\newcommand{\Cb}[1]{\mathrm{C}_{\mathrm{b}}\!\left(#1\right)}
\newcommand{\Cc}[1]{\mathrm{C}_{\mathrm{c}}\!\left(#1\right)}
\newcommand{\C}[1]{\mathrm{C}\!\left(#1\right)}
\newcommand{\refeq}[1]{{\rm (\ref{#1})}}
\newcommand{\proof}{\noindent{\sc Proof.~}}
\newcommand{\qed}{\hfill{\sc Q.E.D.}\medskip}

\newcommand{\Section}[1]{\setcounter{equation}{0}\section{#1}}
\newcommand{\Ma}[1]{\mathscr{M}\!\left(#1\right)}
\newcommand{\Aa}{{\mathscr A}}
\newcommand{\Bb}{{\mathscr B}}
\newcommand{\Ap}{\mathscr{A\!P}}
\newcommand{\AP}{\mathbb{AP}}
\newcommand{\II}{{\mathbb I}}
\newcommand{\JJ}{{\mathbb J}}
\newcommand{\kap}{\kappa}
\newcommand{\kac}[1]{#1_{\text{\tiny\sc kac}}}
\newcommand{\h}{\mathfrak{h}}
\newcommand{\Auhat}{\widehat{A_u}}
\newcommand{\delhat}{\widehat{\del}}
\newcommand{\s}[2]{{#1}_{\scriptscriptstyle\!\:\!#2}}
\newcommand{\sS}[2]{{#1}_{\scriptscriptstyle #2}}

\begin{document}

\title{Quantum Bohr compactification of discrete\\ quantum groups}
\author{Piotr Miko\l{}aj So\l{}tan\thanks{
Partially supported by Komitet Bada\'n Naukowych grants 2 PO3A 040 22
\& 1 PO3A 036 26, the Foundation for Polish Science and Deutsche
Forschungsgemeinschaft.} \\
\small Department of Mathematical Methods in Physics,\\
\small Faculty of Physics, University of Warsaw\\
\small\tt piotr.soltan@fuw.edu.pl}
\maketitle



\abstract{We introduce the analogue of Bohr compactification for discrete
quantum groups on $\cst$-algebra level. The cases of unimodular and general
$\cst$-algebraic discrete quantum groups are treated separately. The passage
from the former case to the latter is done with help of the concept of canonical
Kac quotient for compact quantum groups\footnotemark. Several examples are
presented and passage to the Bohr compactification is investigated.}

\footnotetext{The author is grateful to Stefaan Vaes for introducing him to
this concept.}

\section{Introduction}

The theory of compact quantum groups was initiated with the fundamental papers
\cite{SUq2,pseudogr} of S.L.~Woronowicz. Soon the program of \cite{pseudo} was
carried out for this class of quantum groups and the theory of discrete quantum
groups on the level of $\cst$-algebras was developed in \cite[Section 3]{pw}.
The algebraic counterparts  and generalizations of this theory came in the form
of Van Daele's beautiful theory of duality for multiplier Hopf algebras
with invariant functionals
(\cite{afgd}). It was also Van Daele who found a very appealing description of
discrete quantum groups without reference to their compact duals in
\cite{dqg}. His approach was based on the notion of a multiplier Hopf algebra
which he had introduced in \cite{mha}.

In this paper we intend to use this rich theory to develop a notion of Bohr
compactification of a discrete quantum group on the level of $\cst$-algebras.
It will be a compact quantum group with a universal property mirroring that of
the classical Bohr compactification (\cite{weil,loom,holm}).

The same problem was considered on the algebraic level of multiplier Hopf
$*$-algebras in \cite{cdqg}. The disadvantage of the algebraic approach is the
lack of invariant functionals on the algebraic counterpart of Bohr
compactification. In this paper the use of $\cst$-algebras guarantees the
existence of Haar measure on the Bohr compactification.

Our definition of Bohr compactification follows the link between classical Bohr
compactification and the spectrum of the algebra of almost periodic functions.
As a technical tool we introduce the concept of canonical Kac quotient of a
compact quantum group.

Let us briefly describe the contents of the paper. In Section \ref{defnot} we
collect the basic definitions and introduce notation for the rest of the
paper. We conclude that section with Theorem \ref{unim} giving several
equivalent formulations of the concept of unimodularity for discrete quantum
groups. In the next section we motivate the investigation of Bohr
compactification on $\cst$-algebra level and proceed to define it first in the
unimodular case and then for general $\cst$-algebraic discrete quantum groups.
We formulate the basic properties of quantum Bohr compactification.

In Section \ref{Ex4} we present the first examples of our construction and
relate it to the classical Bohr compactification of discrete groups. We
indicate several corollaries of the universal property of Bohr
compactification and apply them to the relation between the universal
representation of a compact group and its full $\cst$-algebra.

Section \ref{kacsect} is devoted to introduction and study of the canonical Kac
quotient of a compact quantum group. This is an essential tool for the
following considerations. In the next section we give proofs of the results
formulated in Section \ref{QBC} and obtain a useful Corollary \ref{isoKac}.

Section \ref{Ex7} is entirely devoted to examples of our construction. We find
the Bohr compactification of the dual of the quantum $SU(2)$ group of
S.L.~Woronowicz. Then we investigate a possible definition of a discretization
of a compact quantum group and provide an example of the quantum double torus
of Hajac and Masuda whose ``quantum part'' vanishes in the considered procedure.
Finally we consider the family of Bohr compactifications of duals of the
universal compact quantum groups of Van Daele and Wang. We derive its
universal property and prove that all these groups are infinite in an
appropriate sense. We also provide classification of isomorphism classes of
these compact quantum groups showing that the obtained family is not too large.
We conclude that section with a definition of a maximally almost periodic
quantum group. We prove that a $\cst$-algebraic discrete quantum group is
maximally almost periodic if the $\cst$-algebra of functions on its universal
dual is residually finite dimensional. On the other hand we show that if a
$\cst$-algebraic discrete quantum group is maximally almost periodic then it
must be unimodular.

\Section{Definitions and notation}\label{defnot}

In this paper we shall be concerned with the theory of quantum groups. The
main objects of study will be discrete quantum groups. In literature there are
two definitions of discrete quantum groups which are relevant to our
considerations. One comes from the study of compact quantum groups and their
representations (cf.~\cite[Section 3]{pw} and \cite{cqg}) and the other is
formulated in the framework of multiplier Hopf algebras (cf.~\cite{dqg}). In
order to avoid confusion we shall recall them:

\begin{definition}
\noindent
\begin{itemize}
\item[{\rm (1)}] A\/ {\em $\cst$-algebraic discrete quantum group} is a dual of
a compact quantum group.
\item[{\rm (2)}] A\/ {\em discrete quantum group} is a multiplier Hopf
$*$-algebra $(\Aa,\s{\del}{\Aa})$ such that $\Aa$ is an algebraic direct sum of
matrix algebras.
\end{itemize}
\end{definition}

We shall also use two notions of a multiplier algebra. In fact the two notions
coincide, but as the contexts in which they appear are usually different, we
shall make a distinction. For an algebra $\Aa$ with non degenerate product the
{\em algebraic multiplier algebra}\/ will be denoted by $\Ma{\Aa}$. It is the
algebra of double centralizers of $\Aa$ (cf.~\cite{mha}). On the other hand for
a $\cst$-algebra $A$ the {\em $\cst$-algebraic multiplier algebra}\/ of $A$
will be denoted by $\M{A}$. This is the $\cst$-algebra formed by adjointable
maps of the Hilbert $A$-module which is $A$ itself. We remark that all non
degenerate (anti) homomorphisms extend to multiplier algebras (both in the
algebraic and $\cst$-algebraic cases). We shall as a rule denote these
extensions by the same symbol as the original homomorphisms
(\cite[Appendix]{mha}).

Let $(A,\s{\del}{A})$ be a $\cst$-algebraic discrete quantum group. In other
words there exists a compact quantum group $(B,\s{\del}{B})$ such that
$(A,\s{\del}{A})$ is the dual of $(B,\s{\del}{B})$ in the sense of
\cite[Section 3]{pw}. If $\II$ is the set of equivalence classes of all
irreducible unitary representations of $(B,\s{\del}{B})$ then
\begin{equation}\label{dirsum}
A=\bigoplus_{\alpha\in\II} M_{n_\alpha},
\end{equation}
where $n_\alpha$ is the dimension of the representation corresponding to
$\alpha$. This means that $A$ consists of infinite families
$(m_\alpha)_{\alpha\in\II}$ such that $m_\alpha\in M_{n_\alpha}$ and for any
$\eps>0$ there exists a finite subset $F\subset\II$ such that
$\|m_\beta\|\leq\eps$ for all $\beta\in\II\bez F$. The comultiplication
$\s{\del}{A}$ is a coassociative morphism from $A$ to $A\tens A$ (for the
notion of a morphism of $\cst$-algebras we refer the reader to
\cite[Section 0]{unbo}). The multiplier algebra $\M{A}$ consists of those
families $(m_\alpha)_{\alpha\in\II}$ with $m_\alpha\in M_{n_\alpha}$ for which
\[
\sup_{\alpha\in\II}\|m_\alpha\|<\infty.
\]

Let $\Aa$ be the Pedersen ideal of $A$. In view of \refeq{dirsum} it is easy to
see that $\Aa$ is the algebraic direct sum of the same family of matrix
algebras $(M_{n_\alpha})_{\alpha\in\II}$. The product in the algebra  $\Aa$ is
non degenerate and the multiplier algebra $\Ma{\Aa}$ coincides with the set
$\Aff{A}$ of elements affiliated with $A$ in the $\cst$-algebra sense
(\cite[Definition 1.1]{unbo}) and can be identified with the Cartesian product
of the matrix algebras $(M_{n_\alpha})_{\alpha\in\II}$. The situation can be
summed up in the following chain of inclusions:
\[
\Aa\subset A\subset\M{A}\subset\Aff{A}=\Ma{\Aa}.
\]

Let $\s{\del}{\Aa}$ denote the restriction of the map $\s{\del}{A}$ to $\Aa$.
Then $\s{\del}{\Aa}$ is a $*$-homomorphism from $\Aa$ to $\Ma{\Aa}$. Obviously
the range of $\s{\del}{\Aa}$ is contained in $\M{A}$.

Now let $(B,\s{\del}{B})$ be a compact quantum group. Then there exists a
unique dense unital $*$-subalgebra $\Bb\subset B$ such that if $\s{\del}{\Bb}$
denotes the restriction of $\s{\del}{B}$ to $\Bb$ then $(\Bb,\s{\del}{\Bb})$ is
a Hopf $*$-algebra. The algebra $\Bb$ is the linear span of matrix elements of
all finite dimensional unitary representations of $(B,\s{\del}{B})$
(cf.~\cite[Theorem 2.2]{cqg}). We call $(\Bb,\s{\del}{\Bb})$ the Hopf
$*$-algebra {\em associated with $(B,\s{\del}{B})$}.

We include the following proposition for completeness:

\begin{proposition}\label{assoc}
Let $(A,\s{\del}{A})$ be a $\cst$-algebraic discrete quantum group. Let $\Aa$
be the Pedersen ideal of $A$ and let $\s{\del}{\Aa}$ be the restriction of
$\s{\del}{A}$ to $\Aa$. Then the pair $(\Aa,\s{\del}{\Aa})$ is a discrete
quantum group.
\end{proposition}

\proof
Let $(B,\s{\del}{B})$ be a compact quantum group whose dual is
$(A,\s{\del}{A})$. The associated Hopf $*$-algebra $(\Bb,\s{\del}{\Bb})$ of
$(B,\s{\del}{B})$ is a multiplier Hopf $*$-algebra with invariant functionals
as defined in \cite{afgd}. A closer inspection of the beautiful duality theory
presented in that reference shows that $(\Aa,\s{\del}{\Aa})$ is the dual
multiplier Hopf algebra of $(\Bb,\s{\del}{\Bb})$. In particular it is a
multiplier Hopf $*$-algebra. As $\Aa$ is a direct sum of matrix algebras we see
that $(\Aa,\s{\del}{\Aa})$ is a discrete quantum group.
\qed

We shall refer to the discrete quantum group $(\Aa,\s{\del}{\Aa})$, arising
from a $\cst$-algebraic discrete quantum group $(A,\s{\del}{A})$ in the way
described in Proposition \ref{assoc}, as the discrete quantum group
{\em associated with $(A,\s{\del}{A})$}.

Let $(A,\s{\del}{A})$ be a $\cst$-algebraic discrete quantum group and let
$(\Aa,\s{\del}{\Aa})$ be its associated discrete quantum group. The coinverse
$\s{\kap}{A}$ and counit $\s{e}{A}$ restricted to $\Aa$ coincide with the
corresponding maps on the level of $(\Aa,\s{\del}{\Aa})$. The coinverse
$\s{\kap}{\Aa}$ admits a canonical extension to $\Ma{\Aa}=\Aff{A}$. It needs
to be stressed, however, that this map restricted to $\M{A}$ is, in general,
unbounded. Similarly for a compact quantum group $(B,\s{\del}{B})$ and
associated Hopf $*$-algebra $(\Bb,\s{\del}{\Bb})$ the structure maps
$\s{\kap}{\Bb}$ and $\s{e}{\Bb}$ may become unbounded on the level of the
$\cst$-algebra $B$.

Let us end this section with a few remarks about the notion of unimodularity.
We call a discrete quantum group {\em unimodular} if its left and right Haar
measures coincide. The same terminology applies to $\cst$-algebraic discrete
quantum groups. We have the following characterization of unimodularity
(cf.~\cite[Theorem 2.5]{cqg}, \cite[Section 3]{pw}):

\begin{theorem}\label{unim}
Let $(A,\s{\del}{A})$ be a $\cst$-algebraic discrete quantum group and let
$(B,\s{\del}{B})$ be a compact quantum group whose dual is $(A,\s{\del}{A})$.
Then the following statements are equivalent:
\begin{itemize}
\item[{\rm (1)}] $(A,\s{\del}{A})$ is unimodular,
\item[{\rm (2)}] the coinverse $\s{\kap}{A}$ of $(A,\s{\del}{A})$ is bounded,
\item[{\rm (3)}] the coinverse $\s{\kap}{A}$ of $(A,\s{\del}{A})$ is involutive,
\item[{\rm (4)}] either right or left Haar measure of $(A,\s{\del}{A})$ is a
trace,
\item[{\rm (5)}] the coinverse $\s{\kap}{B}$ of $(B,\s{\del}{B})$ is bounded,
\item[{\rm (6)}] the Haar measure of $(B,\s{\del}{B})$ is a trace.
\end{itemize}
\end{theorem}

If a compact quantum group $(B,\s{\del}{B})$ satisfies one of the equivalent
conditions (5) and (6) in Theorem \ref{unim} then we call $(B,\s{\del}{B})$ a
{\em compact quantum group of Kac type.}\/ If the Haar measure of
$(B,\s{\del}{B})$ is faithful then $(B,\s{\del}{B})$ is a
{\em compact Kac $\cst$-algebra} (cf. \cite[D\'efinition 3.2.1]{ev}).

\Section{Quantum Bohr compactification}\label{QBC}

Let $(A,\s{\del}{A})$ be a $\cst$-algebraic discrete quantum group and let
$(\Aa,\s{\del}{\Aa})$ be the associated discrete quantum group. It has been
shown in \cite{cdqg} that the set $\Ap$ of those elements $x$ of $\Ma{\Aa}$ for
which
\[
\s{\del}{\Aa}(x)\in\Ma{\Aa}\algt\Ma{\Aa}
\]
is a Hopf $*$-algebra with the structure maps inherited from $\Ma{\Aa}$. The
elements of $\Ap$ are called {\em almost periodic}\/ elements for
$(\Aa,\s{\del}{\Aa})$. We shall denote this Hopf $*$-algebra by
$(\Ap,\s{\del}{\Aa})$. It has a universal property analogous to the universal
property of Bohr compactification for discrete groups
(cf.~\cite[Theorem 4.5]{cdqg}).

The Hopf $*$-algebra $(\Ap,\s{\del}{\Aa})$ does not in general have positive
invariant functionals. As an example consider the discrete quantum group
$(\Cc{\ZZ},\sS{\del}{\ZZ})$ consisting of the $*$-algebra of finitely supported
functions on $\ZZ$ and the standard comultiplication
$\bigl(\sS{\del}{\ZZ}(f)\bigr)(n_1,n_2)=f(n_1+n_2)$. The element
\[
\mathfrak{n}\colon\ZZ\ni n\longmapsto n\in\CC
\]
is a multiplier of $\Aa$ and belongs to $\Ap$ since
\[
\sS{\del}{\ZZ}(\mathfrak{n})=\mathfrak{n}\tens I+I\tens\mathfrak{n}.
\]
However if $\ph$ were, say, a left invariant positive functional on
$(\Ap,\sS{\del}{\ZZ})$ then $\ph(I)$ would be non zero and
\[
\ph(\mathfrak{n})I=(\id\tens\ph)\sS{\del}{\ZZ}(\mathfrak{n})=
\ph(I)\mathfrak{n}+\ph(\mathfrak{n})I
\]
would yield a contradiction. The problem comes from the unboundedness of
$\mathfrak{n}$. We hope to remedy this problem by restricting to
{\em bounded}\/ almost periodic elements.

\subsection{The unimodular case}\label{unimodular}

Let us first concentrate on the Bohr compactification of a unimodular
$\cst$-algebraic discrete quantum group.

\begin{definition}\label{DefAP}
Let $(A,\s{\del}{A})$ be a unimodular $\cst$-algebraic discrete quantum group
and let $\Ap$ be the algebra of almost periodic elements for the discrete
quantum group associated with $(A,\s{\del}{A})$. Define $\AP$ as the closure in
$\M{A}$ of the set
\[
\Ap\cap\M{A}.
\]
The set $\AP$ is called the\/ {\em space of almost periodic elements for
$(A,\s{\del}{A})$}.
\end{definition}
It is immediate from the definition that $\AP$ is a unital $\cst$-subalgebra of
$\M{A}$ (cf.~\cite[Thm.~4.3]{cdqg}).

\begin{theorem}\label{bialg}
Let $(A,\s{\del}{A})$ be a $\cst$-algebraic discrete quantum group and let
$\AP$ be the $\cst$-subalgebra of\/ $\M{A}$ introduced in Definition
\ref{DefAP}. Then
\begin{equation}\label{teza1}
\s{\del}{A}(\AP)\subset\AP\tens\AP.
\end{equation}
\end{theorem}

\proof
Take $x\in\M{A}\cap\Ap$. Since $\s{\del}{\Aa}$ (extended to $\Ma{\Aa}$)
coincides with the extension of $\s{\del}{A}$ to $\Aff{A}$, we see that
\begin{equation}\label{komn}
\s{\del}{A}(x)=\s{\del}{\Aa}(x).
\end{equation}
The left hand side of \refeq{komn} belongs to $\Ap\algt\Ap$ while the right
hand side is in $\M{A\tens A}$. Let us use these facts:
\begin{equation}\label{sumka}
\s{\del}{\Aa}(x)=\sum_{k=1}^{N}x_k\tens y_k,
\end{equation}
where $x_k,y_k\in\Ap$ for $k=1,\ldots,N$. We can further suppose that the
elements $y_1,\ldots,y_N$ are linearly independent. Consider the family of
functionals $(\xi_l)_{l=1,\ldots,N}$ constructed in the proof of
\cite[Theorem 4.3.1]{cdqg}. They have the property that
\[
\xi_l(y_k)=\delta_{k,l}
\]
and extend to continuous functionals on $A$. In particular we have
\[
(\s{\id}{A}\tens\xi_l)\s{\del}{A}(x)
=(\s{\id}{\Aa}\tens\xi_l)\s{\del}{\Aa}(x)=x_l
\]
for $l=1,\ldots,N$. It follows that the norm of each $x_l$ is finite. Now just
as in the proof of \cite[Theorem 4.3]{cdqg} we can rearrange the sum
\refeq{sumka} in such a way to have a maximal linearly independent subset of
$\{x_1,\ldots,x_N\}$ making up the left leg of $\s{\del}{A}(x)$. Then we can
use the same technique to show that the elements making up the right leg have
finite norms. This way we arrive at the conclusion that $\s{\del}{A}(x)$ can be
written as a sum of elementary tensors on $\Ap\cap\M{A}\subset\AP$. In other
words
\[
\s{\del}{A}(x)\in\AP\algt\AP.
\]
By continuity of $\s{\del}{A}$ we obtain \refeq{teza1}.
\qed

The pair $(\AP,\s{\del}{A})$ is obviously a $\cst$-bialgebra.

\begin{theorem}
Let $(A,\s{\del}{A})$ be a unimodular $\cst$-algebraic discrete quantum
group and let $\AP$ be the algebra introduced in Definition \ref{DefAP}. Then
$(\AP,\s{\del}{A})$ is a compact quantum group.
\end{theorem}

\proof
Recall that $(\Ap,\s{\del}{\Aa})$ is a Hopf $*$-algebra
(\cite[Theorem 4.4]{cdqg}). Therefore the maps
\begin{equation}\label{TT}
\begin{array}
{@{\Ap\algt\Ap\ni a\tens b\longmapsto\;}l@{\smallskip}}
\s{\del}{\Aa}(a)(I\tens b)\in\Ap\algt\Ap,\\
(a\tens I)\s{\del}{\Aa}(b)\in\Ap\algt\Ap
\end{array}
\end{equation}
are bijections (see e.g.~\cite[Prop.~2.1]{mha}). Let $\s{\kap}{\Aa}$ denote
the coinverse be the of the discrete quantum group $(\Aa,\s{\del}{\Aa})$
associated with $(A,\s{\del}{A})$ as well as its extension to
$\Ma{\Aa}=\Aff{A}$. Then coinverse of the Hopf $*$-algebra $(\Ap,\s{\del}{\Aa})$
and the canonical extension of $\s{\kap}{A}$ to $\M{A}$ are restrictions of
$\s{\kap}{\Aa}$. We shall thus use the symbol $\s{\kap}{\Aa}$ also to denote
these two maps.

The inverses of the maps in \refeq{TT} are explicitly given by
\[
\begin{array}
{@{\Ap\algt\Ap\ni a\tens b\longmapsto\;}l@{\smallskip}}
\bigl((\id\tens\s{\kap}{\Aa})\s{\del}{\Aa}(a)\bigr)(I\tens b)\in\Ap\algt\Ap,\\
(a\tens I)\bigl((\s{\kap}{\Aa}\tens\id)\s{\del}{\Aa}(b)\bigr)\in\Ap\algt\Ap,
\end{array}
\]
From the proof of Theorem \ref{bialg} it is clear that if
$x,y\in\M{A}\cap\Ap$ then $\s{\del}{\Aa}(x)(I\tens y)$ and
$(x\tens I)\s{\del}{\Aa}(y)$ are in
$\bigl(\M{A}\cap\Ap\bigr)\algt\bigl(\M{A}\cap\Ap\bigr)$. By Theorem \ref{unim}
the unimodularity of $(A,\s{\del}{A})$ implies that
$\bigl((\id\tens\s{\kap}{\Aa})\s{\del}{\Aa}(x)\bigr)(I\tens y)$ and
$(x\tens I)\bigl((\s{\kap}{\Aa}\tens\id)\s{\del}{\Aa}(y)\bigr)$ also belong to
$\bigl(\M{A}\cap\Ap\bigr)\algt\bigl(\M{A}\cap\Ap\bigr)$, as
$\s{\kap}{\Aa}=\s{\kap}{A}$ is bounded on $\M{A}$. Therefore
$\bigl(\M{A}\cap\Ap\bigr)\algt\bigl(\M{A}\cap\Ap\bigr)$ is in the image of
the maps
\[
\begin{array}
{@{\AP\algt\AP\ni x\tens y\longmapsto\;}l@{\smallskip}}
\s{\del}{\Aa}(x)(I\tens y)=\s{\del}{A}(x)(I\tens y),\\
(x\tens I)\s{\del}{\Aa}(y)=(x\tens I)\s{\del}{A}(y).
\end{array}
\]
Consequently the sets
\[
\begin{array}{l@{\smallskip}}
\bigl\{(x\tens I)\s{\del}{A}(y):\:x,y\in\AP\bigr\},\\
\bigl\{(I\tens x)\s{\del}{A}(y):\:x,y\in\AP\bigr\}
\end{array}
\]
are linearly dense in $\AP\tens\AP$.
\qed

Note that there is no guarantee that $\AP$ is a separable $\cst$-algebra. This
means that we are not exactly in the setting of \cite{cqg}. However, as shown
in e.g.~\cite{mvd}, the separability of the $\cst$-algebra is not essential to
any of the main results of the theory of compact quantum groups.

\begin{definition}\label{dAP}
Let $(A,\s{\del}{A})$ be a unimodular $\cst$-algebraic discrete quantum
group. The compact quantum group $(\AP,\s{\del}{A})$ is called the\/
{\em Bohr compactification} of $(A,\s{\del}{A})$.
\end{definition}

\begin{theorem}\label{matr_elts}
Let $(A,\s{\del}{A})$ be a unimodular $\cst$-algebraic discrete quantum group
and let $(\AP,\s{\del}{A})$ be its Bohr compactification. Then $\AP$ is the
closed linear span of matrix elements of finite dimensional unitary
representations of $(A,\s{\del}{A})$.
\end{theorem}

\proof
Let $u$ be a finite dimensional unitary representation of $(A,\s{\del}{A})$,
i.e.~$u$ is a unitary element of $\M{\K{\CC^N}\tens A}=M_N\bigl(\M{A}\bigr)$
and $(\id\tens\s{\del}{A})u=u_{12}u_{13}$. Let $(u^{kl})_{k,l=1,\ldots,N}$ be
the matrix elements of $u$. Then the familiar formula
\[
\s{\del}{A}(u^{kl})=\sum_{p=1}^Nu^{kp}\tens u^{pl}
\]
shows that all matrix elements of $u$ belong to $\Ap\cap\M{A}\subset\AP$.
Therefore the linear span of matrix elements of all finite dimensional unitary
representations of $(A,\s{\del}{A})$ is contained in $\AP$. On the other Hand
let $v$ be a finite dimensional unitary representation of the compact quantum
group $(\AP,\s{\del}{A})$. Then $v$ is a finite dimensional
unitary representation of $(A,\s{\del}{A})$. It is one of the fundamental
results of the theory of compact quantum groups that matrix elements of finite
dimensional unitary representations of any compact quantum group
$(B,\s{\del}{B})$ span a dense subspace of the $\cst$-algebra $B$. Therefore
the linear span of matrix elements of finite dimensional unitary
representations of $(A,\s{\del}{A})$ is a dense subspace of $\AP$.
\qed

\subsection{The general case}

We shall base our definition of Bohr compactification of a general
$\cst$-algebraic discrete quantum group on Theorem \ref{matr_elts}.

\begin{definition}\label{DefAP_gen}
Let $(A,\s{\del}{A})$ be a $\cst$-algebraic discrete quantum group. We shall
define the space $\AP$ of almost periodic elements for $(A,\s{\del}{A})$ as the
closed linear span of matrix elements of all finite dimensional unitary
representations of $(A,\s{\del}{A})$.
\end{definition}

\begin{theorem}\label{cstar}
Let $(A,\s{\del}{A})$ be a $\cst$-algebraic discrete quantum group and let
$\AP$ be its space of almost periodic elements. Then $\AP$ is a unital
$\cst$-subalgebra of\/ $\M{A}$,
\begin{equation}\label{APinAP}
\s{\del}{A}(\AP)\subset\AP\tens\AP.
\end{equation}
and $(\AP,\s{\del}{A})$ is a compact quantum group.
\end{theorem}

We shall present the proof of this theorem in Section \ref{dowody}.

\begin{definition}\label{DBohr}
Let $(A,\s{\del}{A})$ be a $\cst$-algebraic discrete quantum group and let
$\AP$ be the space of almost periodic elements for $(A,\s{\del}{A})$. The
compact quantum group $(\AP,\s{\del}{A})$ is called the\/ {\em
Bohr compactification} of $(A,\s{\del}{A})$.\end{definition}

Of course, by Theorem \ref{matr_elts} the two definitions of Bohr
compactification coincide for unimodular $\cst$-algebraic discrete quantum
groups.

\begin{remark}\label{algrem}
Let $(A,\s{\del}{A})$ be a $\cst$-algebraic discrete quantum group and let
$(\AP,\s{\del}{A})$ be its Bohr compactification. By definition of $\AP$ the
elements in $\AP$ are limits of multipliers $x$ such that
$\s{\del}{A}(x)\in\M{A}\algt\M{A}$.
\end{remark}

Let us now describe the universal property of Bohr compactification. For a
$\cst$-algebraic discrete quantum group $(A,\s{\del}{A})$ we shall denote by
$\chi$ the inclusion map from $\AP$ into $\M{A}$. It is clear that it is a
unital $*$-homomorphism and hence a morphism of $\cst$-algebras:
$\chi\in\Mor{\AP}{A}$. Moreover it satisfies
$\s{\del}{A}\comp\chi=(\chi\tens\chi)\comp\s{\del}{A}$.

\begin{theorem}\label{univ}
Let $(A,\s{\del}{A})$ be a $\cst$-algebraic discrete quantum group and let
$(\AP,\s{\del}{A})$ be its Bohr compactification. Let $(B,\s{\del}{B})$ be a
compact quantum group. If\/ $\Phi\in\Mor{B}{A}$ is a morphism of quantum
groups, i.e.
\begin{equation}\label{morgr}
\s{\del}{A}\comp\Phi=(\Phi\tens\Phi)\comp\s{\del}{B}.
\end{equation}
Then there exists a unique $\overline{\Phi}\in\Mor{B}{\AP}$ which is a
morphism of compact quantum groups, such that $\Phi=\chi\comp\overline{\Phi}$.
\end{theorem}

We shall give proof of Theorem \ref{univ} in Section \ref{dowody}. It is
clear that the compact quantum group $(\AP,\s{\del}{A})$ is the unique
compact quantum group with the universal property described in Theorem
\ref{univ}.

\Section{First examples}\label{Ex4}

\paragraph{Classical discrete group.}
First we shall see that our definition of Bohr compactification extends the
one known for classical groups.

\begin{proposition}\label{almper}
Let $\Gamma$ be a discrete group and let $(\Ci{\Gamma},\sS{\del}{\Gamma})$ be
the corresponding (unimodular) $\cst$-algebraic discrete quantum group. Let
$(\AP,\sS{\del}{\Gamma})$ be the Bohr compactification of
$(\Ci{\Gamma},\sS{\del}{\Gamma})$ as defined in Definition \ref{dAP}. Then
\begin{itemize}
\item[{\rm (1)}] $\AP$ is isomorphic to the algebra of continuous functions on
the (classical) Bohr compactification of\/ $\Gamma$,
\item[{\rm (2)}] $(\AP,\sS{\del}{\Gamma})$ is the compact quantum group
corresponding to the classical Bohr compactification of\/ $\Gamma$,
\item[{\rm (3)}] $\AP$ can be characterized in the following way:
\begin{equation}\label{almp}
\AP=\bigl\{
f\in\M{\Ci{\Gamma}}:\:\sS{\del}{\Gamma}(f)\in\M{\Ci{\Gamma}}\tens\M{\Ci{\Gamma}}
\bigr\}.
\end{equation}
\end{itemize}
\end{proposition}

\proof
Statements (1) and (2) are consequences of the universal property of
$(\AP,\sS{\del}{\Gamma})$. Indeed: $\AP$ is commutative as it is contained in a
commutative $\cst$-algebra $\M{\Ci{\Gamma}}$. Therefore it corresponds to a
classical compact group. Its universal property for morphisms from compact
quantum groups to $(\Ci{\Gamma},\sS{\del}{\Gamma})$ means in particular that
this compact group has the universal property of the Bohr compactification of
$\Gamma$.

To prove Statement (3) let us first see that the right hand side of
\refeq{almp} is the algebra of almost periodic functions on
$\Gamma$. Indeed: for any locally compact space $X$ and any $\cst$-algebra $C$
the $\cst$-algebra $\M{\Ci{X}}\tens C=\Cb{X}\tens C$ is naturally
isomorphic to the $\cst$-algebra of norm continuous functions from $X$
to $C$ whose image is relatively (norm) compact. Now for any
$f\in\M{\Ci{\Gamma}}=\Cb{\Gamma}$ the property that $\sS{\del}{\Gamma}(f)$ is in
$\Cb{\Gamma}\tens\Cb{\Gamma}$ is equivalent to the statement that the set
\[
\bigl\{f(\cdot\,\gamma):\:\gamma\in\Gamma\bigr\}
\]
is relatively compact in $\Cb{\Gamma}$. This means that $f$ is almost periodic.
It is well known (cf.~\cite[\S 41]{loom}) that the set of all almost periodic
functions on $\Gamma$ is equal to the algebra of continuous functions on the
Bohr compactification of $\Gamma$. This proves the last statement of the
proposition.
\qed

\paragraph{Dual of a compact group.}
Let $G$ be a compact group and $\cst(G)$ its full group $\cst$-algebra. It is
known (cf.~e.g.~\cite[Proposition 3.2]{wn}) that $\cst(G)$ can be made into a
$\cst$-bialgebra by putting $\s{\del}{G}(U_g)=U_g\tens U_g$, where
$G\ni g\mapsto U_g\in\M{\cst(G)}$ is the universal representation of $G$. It
is also well known that this $\cst$-bialgebra is in fact a unimodular
$\cst$-algebraic discrete quantum group.

\begin{theorem}\label{dual}
Let $G$ be a compact group and let $(A,\s{\del}{A})=(\cst(G),\s{\del}{G})$. Let
$(\AP,\s{\del}{A})$ be the Bohr compactification of the (unimodular)
$\cst$-algebraic discrete quantum group $(A,\s{\del}{A})$. Then
$(\AP,\s{\del}{A})$ is isomorphic to
$(\cst(\s{G}{\mathrm{d}}),\del_{G_\mathrm{d}})$, where $\s{G}{\mathrm{d}}$ is
the group $G$ with discrete topology.
\end{theorem}

\proof
Let $(C,\s{\del}{C})$ be the dual of the compact quantum group
$(\AP,\s{\del}{A})$. Since $(A,\s{\del}{A})$ is cocommutative, so is
$(\AP,\s{\del}{A})$ and consequently there exists a discrete group $\Gamma$
such that $(C,\s{\del}{C})=(\cst(\Gamma),\sS{\del}{\Gamma})$ with the usual
compact quantum group structure (cf.~\cite[Theorem 1.7]{pseudogr}). The
morphism $\chi\in\Mor{\AP}{\cst(G)}=\Mor{\cst(\Gamma)}{\cst(G)}$ is a quantum
group homomorphism, i.e.
\[
\s{\del}{G}\comp\chi=(\chi\tens\chi)\comp\sS{\del}{\Gamma}.
\]
In particular for any $\gamma\in\Gamma$ the image of $U_\gamma$ under $\chi$
is a unitary group like element. It follows that $\chi$ determines a group
homomorphism $\xi\colon\Gamma\to G$.

Let $\Theta$ be any discrete group and $(\cst(\Theta),\s{\del}{\Theta})$
the associated compact quantum group (dual of $\Theta$). Then any homomorphism
$\alpha$ from $\Theta$ to $G$ gives a morphism of $\cst$-algebras
$\alpha_*\colon\cst(\Theta)\to\cst(G)$ such that
$\alpha_*(U_\theta)=U_{\alpha(\theta)}$ for all $\theta\in\Theta$. It is easy
to check that this morphism has the property that
\[
\s{\del}{G}\comp\alpha_*=(\alpha_*\tens\alpha_*)\comp\s{\del}{\Theta}.
\]
Therefore there exists a unique
$\overline{\alpha_*}\in\Mor{\cst(\Theta)}{\AP}$ such that
$\alpha_*=\chi\comp\overline{\alpha_*}$
and $\overline{\alpha_*}$ is a compact quantum group morphism.
It follows that
there is a unique homomorphism $\overline{\alpha}\colon\Theta\to\Gamma$
inducing $\overline{\alpha_*}$.

We have thus established that for any discrete group $\Theta$ and any group
homomorphism $\alpha\colon\Theta\to G$ there exists a unique homomorphism
$\overline{\alpha}\colon\Theta\to\Gamma$ such that
$\alpha=\xi\comp\overline{\alpha}$. This means that the pair $(\Gamma,\xi)$
has the universal property with respect to homomorphisms from discrete groups
to $G$. Such a universal pair is clearly unique up to isomorphism. Moreover the
group $\s{G}{\mathrm{d}}$ with the homomorphism
$\id\colon\s{G}{\mathrm{d}}\to G$ has the considered universal property.
Consequently $\Gamma$ is isomorphic to $\s{G}{\mathrm{d}}$ and our theorem
follows.
\qed

Theorem \ref{dual} establishes that the dual of the Bohr compactification of
the dual of a compact group $G$ is the group $G$ with discrete topology. By
Pontryagin duality (cf.~\cite{pseudo}) this means that the Bohr
compactification of $(\cst(G),\s{\del}{G})$ can be equivalently described as
the dual of $\s{G}{\mathrm{d}}$. This fact was well known for abelian groups
(cf.~e.g.~\cite[Sect.~1]{holm}), but without the restriction that $G$ be
compact.

\paragraph{Some corollaries.}
Just as in the algebraic version of Bohr compactification of discrete quantum
groups presented in \cite{cdqg} we can point out an element $\h\in A$ which
does not belong to $\AP$ (cf.~\cite[Prop.~4.7]{cdqg}). It is known that for a
discrete quantum group $(\Aa,\s{\del}{\Aa})$ there exists an element $\h\in\Aa$
such that
\begin{equation}\label{h}
\h a=a\h=\s{e}{\Aa}(a)\h
\end{equation}
for all $a\in\Aa$, where $\s{e}{\Aa}$ is the counit of $(\Aa,\s{\del}{\Aa})$
(\cite[Proposition 3.1]{dqg}).
Now let $(A,\s{\del}{A})$ be a $\cst$-algebraic discrete quantum group and
$(\Aa,\s{\del}{\Aa})$ its associated discrete quantum group. Since
$\Aa\subset A$ the element $\h$ belongs to $A$ and by continuity we have
\refeq{h} for all $a\in A$. It is obvious that $\h$ is unique. What is
important for us is that
\begin{equation}\label{gest}
A\subset\bigl\{(\id\tens\omega)\s{\del}{A}(\h):\:\omega\in A^*\bigr\}^
{\buildrel{\text{\tiny norm}}\over{\text{\tiny closure}}}\hspace{-5ex}.
\end{equation}
Indeed: it is known (cf.~\cite[Section 5]{afgd}) that the set of elements
$(\id\tens\omega_0)\s{\del}{\Aa}(\h)$ with $\omega_0$ a
{\em reduced functional}\/ on $\Aa$ is equal to $\Aa$ (a reduced functional on
$\Aa$ is a functional supported only on finitely many simple summands of
$\Aa$). Since reduced functionals on $\Aa$ extend to continuous functionals on
$A$ and $\Aa$ is dense in $A$, we obtain \refeq{gest}.

\begin{proposition}\label{infi}
Let $(A,\s{\del}{A})$ be a $\cst$-algebraic discrete quantum group with $A$
infinite dimensional and let $(\AP,\s{\del}{A})$ be its Bohr compactification.
Let $\h$ be the element of $A$ with the property that $a\h=\h a=\s{e}{A}(a)\h$
for all $a\in A$. Then $\h$ does not belong to $\AP$.
\end{proposition}

\proof
Suppose that $\h\in\AP$. Then by Remark \ref{algrem} there exists a sequence
$(x_n)_{n\in\NN}$ of multipliers of $A$ such that
$x_n\xrightarrow[n\to\infty]{}\h$ in norm and
\[
\s{\del}{A}(x_n)\in\M{A}\algt\M{A}.
\]
In particular for any $n\in\NN$ the map
\[
A^*\ni\omega\longmapsto(\id\tens\omega)\s{\del}{A}(x_n)\in\M{A}
\]
has finite dimensional range. Consequently the set
\[
\bigl\{(\id\tens\omega)\s{\del}{A}(\h):\:\omega\in A^*,\:\|\omega\|\leq1\bigr\}
\]
is relatively norm compact. However by \refeq{gest} its closure contains the
unit ball of $A$. This shows that the assumption that $\h\in\AP$ was wrong.
\qed

Proposition \ref{infi} has some interesting consequences.

\begin{corollary}
Let $(A,\s{\del}{A})$ be a $\cst$-algebraic discrete quantum
group with $A$ infinite dimensional. Then $(\M{A},\s{\del}{A})$ is not a
compact quantum group.
\end{corollary}

\proof
This is a consequence of the universal property of $(\AP,\s{\del}{A})$
described in Theorem \ref{univ} and of Proposition \ref{infi}. We know that
$\AP$ does not contain $A$, so in particular, it is not equal to $\M{A}$. If
$(\M{A},\s{\del}{A})$ were a compact quantum group then the identity map
$\M{A}\to\M{A}$ would have to factor through the inclusion of $\AP$ into
$\M{A}$. This is not possible.
\qed

For another corollary of Proposition \ref{infi} let us return to the example
discussed in the previous paragraph. Let $G$ be a compact group, $\cst(G)$ its
full group $\cst$-algebra and let
\begin{equation}\label{odwz}
G\ni g\longmapsto U_g\in\M{\cst(G)}.
\end{equation}
be the universal representation of $G$. It is known that $\cst(G)$ is generated
by the elements $\{U_g:\:g\in G\}$ in the sense of \cite{gen}
(see \cite[Theorem 3.1]{wn}).\footnote{In fact the map \refeq{odwz} is a
unitary element of $\M{\Ci{G}\tens\cst(G)}$ which constitutes a quantum family
of elements affiliated with $\cst(G)$ generating this $\cst$-algebra (see
\cite[Section 4]{gen}).}
The relation between the generating elements and the algebra $\cst(G)$ is,
however, rather complicated.

\begin{corollary}\sloppy
Let $G$ be an infinite compact group. Then the $\cst$-subalgebra of
$\M{\cst(G)}$ generated by the universal representation of $G$ does
not contain $\cst(G)$.
\end{corollary}

\proof
Let $(A,\s{\del}{A})=(\cst(G),\s{\del}{G})$ be the $\cst$-algebraic discrete
quantum group obtained from $G$ (the dual of $G$). Also let $(\AP,\s{\del}{A})$
be the Bohr compactification of this $\cst$-algebraic discrete quantum group.
It is immediate that all elements $(U_g)_{g\in G}$ are contained in $\AP$ (they
are matrix elements of one dimensional representations of
$(\cst(G),\s{\del}{G})$). Therefore the $\cst$-subalgebra of $\M{\cst(G)}$
generated by these elements is a subalgebra of $\AP$. By Proposition \ref{infi}
it does not contain $A$.
\qed

\Section{Canonical Kac quotient of a compact quantum group}\label{kacsect}

In this section we shall discuss the concept of the canonical Kac qoutient of a
compact quantum group. This notion is due to Stefaan Vaes \cite{va}.

Let $(B,\s{\del}{B})$ be a compact quantum group. Let $J$ be the (closed two
sided) ideal of $B$ defined as the intersection of left kernels of all tracial
states on $B$:
\begin{equation}\label{J}
J=\bigl\{b\in B:\:\tau(b^*b)=0\text{ for any tracial state $\tau$ on }B\bigr\}
\end{equation}
(in case there are no tracial states we set $J=B$). Let $\kac{B}$ be the
quotient $B/J$ and let $\pi$ be the corresponding quotient map.

\begin{proposition}\label{KacProp}
Let $(B,\s{\del}{B})$ be a compact quantum group and let $\kac{B}$ be the
quotient of $B$ by the ideal \refeq{J} with $\pi$ the quotient map. Then
\begin{itemize}
\item[{\rm (1)}] The equation
\begin{equation}\label{kacdel}
\kac{\del}\bigl(\pi(b)\bigr)=(\pi\tens\pi)\s{\del}{B}(b)
\end{equation}
defines a comultiplication on $\kac{B}$; with this comultiplication
$(\kac{B},\kac{\del})$ becomes a compact quantum group;
\item[{\rm (2)}] $(\kac{B},\kac{\del})$ is a compact quantum group of Kac
type.
\end{itemize}
\end{proposition}

\proof
First let us remark that $\kac{B}$ is evidently isomorphic to the image of $B$
in the representation which is the direct sum of GNS representations of
$B$ for all tracial states. In particular $\kac{B}$ has a faithful family of
tracial states (and if it is separable then it possesses a faithful tracial
state). Moreover $\ker{(\pi\tens\pi)}$ consists of all $x\in B\tens B$ such
that $\bigl(\tau_1\tens\tau_2\bigr)(x^*x)=0$ for any tracial states
$\tau_1,\tau_2$ of $B$.

{\sc Ad (1).} Let $b\in B$ be such that $\pi(b)=0$. Then for any two tracial
states $\tau_1,\tau_2$ on $B$ we have
\[
\bigl(\tau_1\tens\tau_2\bigr)\bigl(\s{\del}{B}(b)^*\s{\del}{B}(b)\bigr)=
\bigl(\tau_1*\tau_2\bigr)(b^*b)
\]
which is equal to $0$, since a convolution of two traces is a trace (and
$b\in J$). Therefore \refeq{kacdel} defines a map
$\kac{B}\to\kac{B}\tens\kac{B}$. It is now straightforward to check that
$\kac{\del}$ is coassociative. Moreover since for any $a,b\in B$
\[
\begin{array}{r@{\;=\;}l@{\smallskip}}
(\pi\tens\pi)\bigl(\s{\del}{B}(a)(I\tens b)\bigr)
&\kac{\del}\bigl(\pi(a)\bigr)\bigl(I\tens\pi(b)\bigr),\\
(\pi\tens\pi)\bigl((a\tens I)\s{\del}{B}(b)\bigr)
&\bigl(I\tens\pi(a)\bigr)\kac{\del}\bigl(\pi(b)\bigr)
\end{array}
\]
and $\pi\tens\pi$ is surjective, we see that $(\kac{B},\kac{\del})$ is a
compact quantum group.

{\sc Ad (2).} We shall show that the Haar measure of $(\kac{B},\kac{\del})$ is
a trace and the conclusion will follow from Theorem \ref{unim}. To that end we
shall repeat the procedure of constructing the Haar measure described in
\cite[Section 4]{mvd} and use a slight modification of the argument in
\cite[Lemma 3.1]{cqg} (for separable $\kac{B}$ it is enough to inspect the
proof of \cite[Theorem 2.3]{cqg} or \cite[Theorem 4.2]{pseudogr}).

We shall use the following generalization of \cite[Lemma 3.1]{cqg}:
\medskip

{\em Let $(D,\s{\del}{D})$ be a compact quantum group and let
$(\rho_\iota)_{\iota\in\JJ}$ be a faithful family of states of $D$. Let $h$ be
a state of $D$ such that
\[
h*\rho_\iota=\rho_\iota*h=h
\]
for all $\iota\in\JJ$. Then $h$ is the Haar measure of $(D,\s{\del}{D})$.}
\medskip

The proof of the above fact is the same as \cite[Lemma 3.1]{cqg}. In the
original formulation the family $(\rho_\iota)_{\iota\in\JJ}$ consisted of a
single element.

Now we can repeat the argument of A.~Maes and A.~Van Daele from \cite{mvd}.
For any tracial state $\tau$ of $\kac{B}$ there exists a tracial state $h$ such
that $h*\tau=\tau*h=h$. This is \cite[Lemma 4.2]{mvd} combined with the fact
that a convolution of traces is a trace.

Let $\tau$ be any tracial state of $\kac{B}$ and let $h$ be a tracial state on
$\kac{B}$ such that $h*\tau=\tau*h=h$. Then if $\omega$ is a positive
functional on $\kac{B}$ such that $\omega\leq\tau$ then by \cite[Lemma 4.3]{mvd}
we have that $h*\omega=\omega*h=\omega(I)h$.

For any tracial positive functional $\rho$ we set
\[
K_\rho=\bigl\{h:\: h\text{ is a tracial state of }\kac{B}
\text{ such that }h*\rho=\rho*h=\rho(I)h\bigr\}.
\]
Then $K_\rho$ is non empty and weakly compact. As in \cite{mvd} we have
$K_{\rho_1+\rho_2}\subset K_{\rho_1}\cap K_{\rho_2}$ and the family of all
these sets has non empty intersection. Let $\kac{h}$ be an element of this
intersection. Then $\kac{h}$ satisfies
\[
\kac{h}*\tau=\tau*\kac{h}=\kac{h}
\]
for all tracial states $\tau$. Since $\kac{B}$ has a faithful family of tracial
states, it follows that $\kac{h}$ is the Haar measure of $(\kac{B},\kac{\del})$.
By construction $\kac{h}$ is a trace.
\qed

We shall refer to the compact quantum group of Kac type $(\kac{B},\kac{\del})$
as the {\em canonical Kac quotient}\/ of the compact quantum group
$(B,\s{\del}{B})$. This terminology is not fully consistent with our approach
to quantum groups because $(\kac{B},\kac{\del})$ corresponds to a
{\em subgroup}\/ of $(B,\s{\del}{B})$. Nevertheless we have decided to use it
because the feature of $(\kac{B},\kac{\del})$ we want to use is related to
properties of $\kac{B}$ as a $\cst$-algebra. More precisely we shall use the
easy fact that any finite dimensional representation of the $\cst$-algebra $B$
factors through the map $\pi\colon B\to\kac{B}$.\footnote{Of course this is
even true for any representation generating a finite von Neumann algebra.}
Notice, however, that the map $\pi\in\Mor{B}{\kac{B}}$ is a quantum group
morphism.

\begin{remark}\label{KacRem}
In the proof of statement (2) of Proposition \ref{KacProp} we have shown that
if a compact quantum group $(B,\s{\del}{B})$ has the property that $B$
possesses a faithful family of tracial states then $(B,\s{\del}{B})$ is its
own canonical Kac qoutient. In particular if $(B,\s{\del}{B})$ is a compact Kac
$\cst$-algebra then its Haar measure is faithful and tracial. Therefore the
canonical Kac quotient of $(B,\s{\del}{B})$ is equal to $(B,\s{\del}{B})$.
\end{remark}

A particular example of the situation described in Remark \ref{KacRem} is the
following: recall that a $\cst$-algebra $B$ is called {\em residually finite
dimensional}\/ if it possesses a separating family of finite dimensional
representations. Now let $(B,\s{\del}{B})$ be a compact quantum group and let
$\pi$ be the quotient map from $B$ onto $\kac{B}$. If $B$ is residually finite
dimensional then it possesses a faithful family of tracial states and
consequently $\pi$ is an isomorphism. As a corollary we get:

\begin{corollary}\label{RDFkac}
Let $(B,\s{\del}{B})$ be a compact quantum group with $B$ a residually finite
dimensional $\cst$-algebra. Then $(B,\s{\del}{B})$ is of Kac type.
\end{corollary}

\Section{Proofs}\label{dowody}

\paragraph{Proof of Theorem \ref{cstar}.}
We shall be very explicit in the proof and use solely the most basic
reference \cite{pw}.

Let $(A,\s{\del}{A})$ be a $\cst$-algebraic discrete quantum group and let
$(B,\s{\del}{B})$ be its universal dual. Let $(\kac{B},\kac{\del})$ be the
canonical Kac quotient of $(B,\s{\del}{B})$ and let $(C,\s{\del}{C})$ be the
dual of $(\kac{B},\kac{\del})$. Then $(C,\s{\del}{C})$ is a unimodular
$\cst$-algebraic discrete quantum group. Denote by $\pi$ the quotient map
$B\to\kac{B}$. Let $u\in\M{A\tens B}$ and $v\in\M{C\tens\kac{B}}$ be the
universal bicharacters (\cite[Section 2]{pw}) for the dualities
between $(B,\s{\del}{B})$ and $(A,\s{\del}{A})$ and between
$(\kac{B},\kac{\del})$ and $(C,\s{\del}{C})$ respectively. We can summarize all
this information in the following diagram:
\[
\xymatrix{(A,\s{\del}{A})\ar@{<->}[rr]^-{\text{\tiny Duality}}_-{u}
&&(B,\s{\del}{B})\ar[d]^{\pi}\\
&&(\kac{B},\kac{\del})\ar@{<->}[rr]^-{\text{\tiny
Duality}}_-{v}&&(C,\s{\del}{C})}
\]
Put $\kac{u}=(\id\tens\pi)u\in\M{A\tens\kac{B}}$. Then using
\cite[Theorem 2.1]{pw} we easily show that there exists a unique
$\Psi\in\Mor{C}{A}$ such that $\kac{u}=(\Psi\tens\id)v$. The morphism $\Psi$
is sometimes called the {\em dual}\/ of $\pi$.

Let us examine the morphism $\Psi$ in more detail. Notice that since $\pi$ is a
quantum group morphism, from any irreducible unitary representation $w$ of
$(B,\s{\del}{B})$ we can make a finite dimensional unitary representation of
$(\kac{B},\kac{\del})$ by applying $\pi$ to each matrix entry of $w$.
We can decompose it into irreducible blocks and then map the corresponding
matrix algebras contained in $C$ into the matrix algebra generated by the left
leg of $w$ inside $A$. As $\pi$ is surjective, we see that any irredicible
representation of $(\kac{B},\kac{\del})$ will appear in decomposition of some
$w$. This way
we obtain an injective mapping from $C$ to $\M{A}$ (cf.~\cite[Section 3]{pw} and
Section \ref{defnot}). A moment of reflection shows that this mapping applied
to the left leg of $v$ gives exactly $\kac{u}$. Since $\Psi$ has the same
property and this property determines $\Psi$ uniquely, it follows that $\Psi$
is injective.

It is also easy to see that $\Psi$ is a quantum group morphism. Indeed the maps
$\Phi_1=\s{\del}{A}\comp\Psi$ and $\Phi_2=(\Psi\tens\Psi)\comp\s{\del}{C}$ have
the property that
\[
(\Phi_k\tens\id)v=(\kac{u})_{23}(\kac{u})_{13}
\]
for $k=1,2$ and the uniqueness part of \cite[Theorem 2.1]{pw} shows that they
are equal.

We have established that
\[
(\id\tens\pi)u=\kac{u}=(\Psi\tens\id)v.
\]
Let $w\in M_N\bigl(\M{A}\bigr)$ be a finite dimensional unitary representation
of $(A,\s{\del}{A})$. Then there exists a unique representation
$\xi\in\Rep{B}{\CC^N}$ such that $w=(\id\tens\xi)u^*$
(cf.~\cite[Section 3]{pw}). Also there is a unique
$\kac{\xi}\in\Rep{\kac{B}}{\CC^N}$ such that $\kac{\xi}\comp\pi=\xi$. This
representation, in turn, gives rise to $\widetilde{w}=(\id\tens\kac{\xi})v^*$
which is a finite dimensional representation of $(C,\s{\del}{C})$. We have
\[
\begin{array}{r@{\;=\;}l@{\smallskip}}
(\Psi\tens\id)\widetilde{w}
&(\Psi\tens\kac{\xi})v^*\\
&(\id\tens\kac{\xi})(\Psi\tens\id)v^*\\
&(\id\tens\kac{\xi})(\kac{u})^*\\
&(\id\tens\kac{\xi})(\id\tens\pi)u^*\\
&(\id\tens\xi)u^*=w.
\end{array}
\]
It follows that any matrix element of $w$ is an image under $\Psi$ of a matrix
element of a finite dimensional unitary representation of $(C,\s{\del}{C})$.
Therefore the closed linear span of matrix elements of all finite dimensional
unitary representations of $(A,\s{\del}{A})$ coincides with the image under
$\Psi$ of the closed linear span of matrix elements of all finite dimensional
unitary representations of $(C,\s{\del}{C})$. By Theorem \ref{matr_elts}, this
last set is a unital $\cst$-subalgebra of $\M{C}$ (since $(C,\s{\del}{C})$ is
unimodular). It follows that $\AP$ is a unital $\cst$-subalgebra of $\M{A}$.

Formula \refeq{APinAP} holds by definition of the $\cst$-algebra $\AP$.
Moreover, since $\Psi$ is a quantum group morphism, we see that the density
conditions for $(\AP,\s{\del}{A})$ follow from those for the Bohr
compactification of $(C,\s{\del}{C})$. This concludes the proof of Theorem
\ref{cstar}.
\qed

\paragraph{Proof of Theorem \ref{univ}.}
Let $(B,\s{\del}{B})$ be a compact quantum group and let $\Phi\in\Mor{B}{A}$
be a quantum group morphism. For any finite dimensional unitary representation
$u\in M_N(B)$ of $(B,\s{\del}{B})$ the matrix
$(\id\tens\Phi)u\in M_N\bigl(\M{A}\bigr)$ is a finite dimensional unitary
representation of $(A,\s{\del}{A})$ (by \refeq{morgr}). Since the span of
matrix elements of finite dimensional unitary representations of
$(B,\s{\del}{B})$ is dense in $B$, we see that $\Phi(B)\subset\AP$. Let
$\overline{\Phi}$ be the map $\Phi$ considered as a morphism from $B$ to $\AP$.
It is straightforward that $\overline{\Phi}$ is a compact quantum group
morphism and that $\Phi=\chi\comp\overline{\Phi}$. The uniqueness of
$\overline{\Phi}$ follows from injectivity of $\chi$.
\qed

\begin{corollary}\label{isoKac}
Let $(A,\s{\del}{A})$ be a $\cst$-algebraic discrete quantum group and let
$(\AP,\s{\del}{A})$ be its Bohr compactification. Then $(\AP,\s{\del}{A})$ is a
compact quantum group of Kac type. Moreover $(\AP,\s{\del}{A})$ is isomorphic
to the Bohr compactification of the dual of the canonical Kac quotient of the
universal dual of $(A,\s{\del}{A})$.
\end{corollary}

\proof
First let us notice that if $(C,\s{\del}{C})$ is a unimodular $\cst$-algebraic
discrete quantum group then its Bohr compactification is a compact quantum
group of Kac type. This is because the coinverse of the Bohr compactification is
a restriction of the coinverse of $(C,\s{\del}{C})$ extended to $\M{C}$. This
extension is still involutive and hence bounded. Thus by Theorem \ref{unim}
the Bohr compactification of $(C,\s{\del}{C})$ is a compact quantum group of Kac
type.

As we have seen in the proof of Theorem \ref{cstar} the compact quantum group
$(\AP,\s{\del}{A})$ is equal to the image under an injective compact quantum
group morphism of the Bohr compactification of the dual of the canonical Kac
quotient $(\kac{B},\kac{\del})$ of the universal dual $(B,\s{\del}{B})$ of
$(A,\s{\del}{A})$. In particular, since the dual of $(\kac{B},\kac{\del})$ is
unimodular, $(\AP,\s{\del}{A})$ is isomorphic to the Bohr compactification of a
unimodular $\cst$-algebraic discrete quantum group and hence a compact quantum
group of Kac type.
\qed

\begin{remark}
Let $(A,\s{\del}{A})$ be a $\cst$-algebraic discrete quantum group. Then the
$\cst$-algebra $\AP$ of almost periodic elements for $(A,\s{\del}{A})$ can be
defined directly as the closure in $\M{A}$ of the set of those multipliers $x$
which lie in $\Ap$ (cf.~Subsection \ref{unimodular}) and such that both
$\s{\kap}{A}(x)$ and $\s{\kap}{A}(x^*)$ are in $\M{A}$. In case
$(A,\s{\del}{A})$ is unimodular this definition reduces to Definition
\ref{DefAP} by Theorem \ref{unim}.
\end{remark}

\Section{Further examples}\label{Ex7}

\paragraph{The dual of quantum $SU(2)$.} Let $0<q<1$. Then the canonical Kac
quotient of the compact quantum group $S_qU(2)$ is isomorphic to the compact
quantum group $(\C{\TT},\sS{\del}{\TT})$. Indeed: recall (\cite{SUq2}) that the
algebra $B$ of continuous functions on $S_qU(2)$ is generated by two elements
$\alpha$ and $\gamma$ with the following commutation relations:
\[
\begin{array}{r@{\;}c@{\;}l@{\smallskip}}
\alpha\gamma&=&q\gamma\alpha,\\
\alpha\gamma^*&=&q\gamma^*\alpha,\\
\gamma^*\gamma&=&\gamma\gamma^*,\\
\alpha^*\alpha&+&\gamma^*\gamma=I,\\
\alpha\alpha^*&+&q^2\gamma^*\gamma=I.
\end{array}
\]
It is immediate from these relations that any tracial state of $B$ must vanish
on $\gamma^*\gamma$. It follows that $\kac{B}$ is the universal
$\cst$-algebra generated by $\pi(\alpha)$ with the relation that $\pi(\alpha)$
is unitary. Also $S_qU(2)$ is the universal dual of
$\widehat{S_qU}(2)$. Therefore using Corollary \ref{isoKac} we see
that the Bohr compactification of $\widehat{S_qU}(2)$ is isomorphic to the
classical compact group arising as the Bohr compactification of $\ZZ$.

\paragraph{Discretization of the double torus.} The quantum double torus was
constructed by P.M.~Hajac and T.~Masuda in \cite{hm}. It is a compact quantum
group arising as the subgroup of the quantum general linear group
$\mathrm{GL}_{\overline{q},q}(2,\CC)$ with $|q|=1$ ``preserving the
noncommutative torus'' and it extends the action of the classical torus on the
noncommutative torus (see \cite[Section 1]{hm} for details).

Let $q$ be a complex number of absolute value 1. The quantum double torus is a
compact quantum group $(B,\s{\del}{B})$ such that
$B=\C{\TT^2}\oplus\mathcal{A}_q$, where $\mathcal{A}_q$ denotes the rotation
algebra (non commutative torus) generated by two unitary elements $b$ and $c$
satisfying $bc=q^2cb$. Let $a$ and $d$ be the standard generators of the
algebra of continuous functions on the classical torus. Then $a,b,c$ and $d$
are elements of $B$ and they generate this $\cst$-algebra. The
comultiplication $\s{\del}{B}$ acts on generators in the following way:
\[
\begin{array}{r@{\;=\;}l@{\smallskip}}
\s{\del}{B}(a)&a\tens a+b\tens c,\\
\s{\del}{B}(b)&a\tens b+b\tens d,\\
\s{\del}{B}(c)&c\tens a+d\tens c,\\
\s{\del}{B}(d)&c\tens b+d\tens d.
\end{array}
\]
The pair $(B,\s{\del}{B})$ is a compact quantum group and it is called the
{\em quantum double torus}. Since $B$ possesses a faithful trace, the quantum
double torus is a compact Kac $\cst$-algebra. In fact $(B,\s{\del}{B})$
is a compact {\em matrix}\/ quantum group with fundamental representation
\[
\begin{pmatrix} a&b\\ c&d\end{pmatrix}.
\]

Theorem \ref{dual} suggests that one could define the discretization of the
quantum space underlying $(B,\s{\del}{B})$ as the $\cst$-algebra of the dual
quantum group of the Bohr compactification of the dual of $(B,\s{\del}{B})$.
However this construction of a ``quantum space with discrete topology'' may
result in a great loss of information. This is indeed the case with the
quantum double torus.

As $(B,\s{\del}{B})$ is a compact Kac $\cst$-algebra, it is isomorphic to
its canonical Kac quotient, by Remark \ref{KacRem}. Let $(A,\s{\del}{A})$ be
the dual of $(B,\s{\del}{B})$
(described in detail in \cite[Section 3]{hm}) and let $(\AP,\s{\del}{A})$ be
its Bohr compactification. Then $\AP$ is the direct sum of algebras
$\B{H_\rho}$, where $\rho$ runs over all pairwise inequivalent irreducible
finite dimensional representations of $B$ and $H_\rho$ is the carrier space of
$\rho$.

If $q=\exp(2\pi i\theta)$ with $\theta$ irrational then $\mathcal{A}_q$ has no
finite dimensional representations and consequently the discretization of the
quantum double torus is then simply a classical (single) 2-torus with
discrete topology. For rational $\theta$ the algebra $\mathcal{A}_q$ is
isomorphic to $M_N\tens\C{\TT}$ and consequently the discretization of the
quantum double torus is the disjoint union (direct sum on the level of
algebras of functions) of the discretized classical 2-torus and a matrix
algebra over discretized circle.

In both cases described above the discretization of the quantum space
underlying $(B,\s{\del}{B})$ turns out to be commutative or almost
commutative. The important information carried by the rotation algebra (for
irrational $\theta$) is lost.

\paragraph{Duals of universal compact quantum groups.}
The universal quantum groups were introduced by A.~Van Daele and S.~Wang in
\cite{vawa} with a special subclass described earlier in \cite{free}. It is a
family of compact matrix quantum groups $(A_u(Q),\s{\del}{Q})$ parametrized by
nonsingular complex matrices $Q$. For a given matrix $Q\in M_m$, the
$\cst$-algebra $A_u(Q)$ is generated by entries of a matrix
$u\in M_m\bigl(A_u(Q)\bigr)=M_m\tens A_u(Q)$ subject to relations
\begin{equation}\label{unirel}
\begin{array}{r@{\:=I_m\tens I=\;}l@{\smallskip}}
u^*u&uu^*,\\
u^\top(Q\tens I)\overline{u}(Q^{-1}\tens I)
&(Q\tens I)\overline{u}(Q^{-1}\tens I)u^\top
\end{array}
\end{equation}
where the $(k,l)$-entry of $\overline{u}$ is $u_{k,l}^*$ and the $(k,l)$-entry
of $u^\top$ is $u_{l,k}$. The universal property of this family of compact
quantum groups is that any compact matrix quantum group is a {\em quantum
subgroup}\/ of one of the universal compact quantum groups (see
\cite[Theorem 1.3 (2)]{vawa}).

In this paragraph we shall be concerned with the family of Bohr
compactifications of the duals of universal quantum groups. For a nonsingular
$m\times m$ complex matrix $Q$ let $(\Auhat(Q),\s{\delhat}{Q})$ be the dual of
$(A_u(Q),\s{\del}{Q})$ and let $(\s{\AP}{Q},\s{\del}{Q})$ be the Bohr
compactification of $(\Auhat(Q),\s{\delhat}{Q})$. An immediate consequence of
theuniversal property of the family of the compact quantum groups
$\bigl\{(A_u(Q),\s{\del}{Q})\bigr\}$ is that for any $\cst$-algebraic discrete
quantum group $(A,\s{\del}{A})$ whose dual is a compact {\em matrix}\/ quantum
group the $\cst$-algebra $A$ embeds into $\M{\Auhat(Q)}$ for some matrix $Q$
with a comultiplication preserving morphism (cf.~\cite[Theorem 1.3]{vawa}).

Let $(F,\s{\del}{F})$ be a finite quantum group (compact quantum group with
$F$ finite dimensional). Then the dual of
$(F,\s{\del}{F})$ is also a finite quantum group. In particular it is a
$\cst$-algebraic discrete quantum group whose dual is a compact matrix quantum
group. Therefore $(F,\s{\del}{F})$ embeds into $(\Auhat(Q),\s{\delhat}{Q})$ for
some matrix $Q$. But at the same time it is a compact quantum group, so the
universal property of Bohr compactification implies that $(F,\s{\del}{F})$
embeds into $(\s{\AP}{Q},\s{\del}{Q})$. This means that the family of compact
quantum groups of the form $(\s{\AP}{Q},\s{\del}{Q})$ has the following
universal property:

\begin{proposition}\label{QQ}
For each $Q\in\mathrm{GL}(m,\CC)$ let $(\s{\AP}{Q},\s{\del}{Q})$ be the Bohr
compactification of the dual $(\Auhat(Q),\s{\delhat}{Q})$ of the universal
quantum group $(A_u(Q),\s{\del}{Q})$.
Let $(F,\s{\del}{F})$ be a finite quantum group. Then there exist an
$m\in\NN$ and $Q\in\mathrm{GL}(m,\CC)$ such that $(F,\s{\del}{F})$ is a
quotient of $(\s{\AP}{Q},\s{\del}{Q})$ in the sense that $F$ embeds into
$\s{\AP}{Q}$ with a compact quantum group morphism.
\end{proposition}

The family $\bigl\{(\s{\AP}{Q},\s{\del}{Q})\bigr\}_{\scriptscriptstyle
Q\in\mathrm{GL}(m,\CC),\:m\in\NN}$ described in
Proposition \ref{QQ} can be called a family of {\em profinite quantum
groups.}\footnote{This terminology was suggested by Shuzhou Wang.}

Clearly it is not difficult to obtain a family of compact quantum groups
satisfying the universal property described in Proposition \ref{QQ}. As an
example one can take the family of all finite quantum groups. However we shall
show that the quantum groups $\bigl\{(\s{\AP}{Q},\s{\del}{Q})\bigr\}_{
\scriptscriptstyle Q\in\mathrm{GL}(m,\CC),\:m\in\NN}$ are all infinite
(i.e.~the algebras $\bigl\{\s{\AP}{Q}\bigr\}$ are infinite dimensional). In
fact the algebras $\s{\AP}{Q}$ are all non separable. We shall also prove that
there are many isomorphic quantum groups in the family
$\bigl\{(\s{\AP}{Q},\s{\del}{Q})\bigr\}_{\scriptscriptstyle
Q\in\mathrm{GL}(m,\CC),\:m\in\NN}$
and classify the isomorphism classes. For this we shall need some easy facts
about quantum subgroups of compact quantum groups.

Consider a compact quantum group $(B^1,\del^1)$. A {\em quantum subgroup}\/ of
$(B^1,\del^1)$ is a compact quantum group $(B^2,\del^2)$ with a
$\cst$-epimorphism $\Theta\colon B^1\to B^2$ such that
$\del^2\comp\Theta=(\Theta\tens\Theta)\comp\del^1$. If $\tau$ is any tracial
state on $B^2$ then $\tau\comp\Theta$ is a tracial state on $B^1$.
Consequently if $b\in B^1$ has the property that $b^*b$ is sent to $0$ by any
tracial state of $B^1$ then $\Theta(b)^*\Theta(b)$ is mapped to zero by any
tracial state of $B^2$. Let $(\kac{B^1},\kac{\del^1})$ and
$(\kac{B^2},\kac{\del^2})$  be the canonical Kac quotients of $(B^1,\del^1)$ and
$(B^2,\del^2)$ respectively and let $\pi_k\colon B^k\to\kac{B^k}$ ($k=1,2$) be
the corresponding quotient maps. It follows that there exists a unique
$\cst$-epimorphism $\widetilde{\Theta}\colon\kac{B^1}\to\kac{B^2}$ such that
the diagram
\[
\xymatrix{B^1\ar[rr]^-{\Theta}\ar[d]^{\pi_1}&&
B^2\ar[d]^{\pi_2}\\
\kac{B^1}\ar[rr]^-{\widetilde{\Theta}}&&\kac{B^2}}
\]
is commutative. It is easy to see that $\widetilde{\Theta}$ is a morphism of
quantum groups.

\begin{lemma}\label{klasy}
For $k=1,2$ let $(B^k,\del^k)$ be a compact quantum groups with canonical Kac
quotients $(\kac{B^k},\kac{\del^k})$ and let $\pi_k\colon B^k\to\kac{B^k}$ be
the quotient maps. Assume that $(B^2,\del^2)$ is a quantum subgroup of
$(B^1,\del^1)$ with the $\cst$-epimorphism $\Theta\colon B^1\to B^2$ and that
$\ker{\Theta}\subset\ker{\pi_1}$. Then the induced morphism
$\widetilde{\Theta}\colon\kac{B^1}\to\kac{B^2}$ is an isomorphism of compact
quantum groups.
\end{lemma}

\proof
Let $\tau$ be a tracial state on $B^1$ and let $(H,\rho,\Omega)$ be the
associated GNS triple. The representation $\rho$ factors through $\Theta$,
i.e.~there is a representation $\widetilde{\rho}$ of $B^2$ in the Hilbert
space $H$ such that $\rho=\widetilde{\rho}\comp\Theta$. This is because
$\ker{\Theta}\subset\ker{\pi_1}\subset\ker{\rho}$.  It follows that
$\tau=\widetilde{\tau}\comp\Theta$ where $\widetilde{\tau}$ is the tracial
state of $B^2$ given by
\[
\widetilde{\tau}(y)=\its{\Omega}{\widetilde{\rho}(y)}{\Omega}
\]
(it is easy to check that this is a trace).  It follows that any trace on
$B^1$ is a composition of $\Theta$ and a trace on $B^2$.

Suppose that we have an element $z=\pi_1(x)$ such that
$\widetilde{\Theta}(z)=0$. Then $\Theta(x)$ lies in the kernel of $\pi_2$.
This means that any tracial state of $B^2$ vanishes on $\Theta(x^*x)$. In view
of the reasoning above all tracial states of $B^1$ vanish on $x^*x$. In other
words $x\in\ker{\pi_1}$ and consequently $z=0$. This shows that the
$\cst$-epimorphism $\widetilde{\Theta}$ is injective.
\qed

We shall use Lemma \ref{klasy} to describe the isomorphism classes of the
family of Bohr compactifications of the duals of universal quantum groups.

\sloppy
Recall (\cite[Section 1]{banica} and \cite{iso}) that we have quantum group
isomorphisms between $(A_u(Q),\s{\del}{Q})$ and
$(A_u(|Q|),\s{\del}{|Q|})$ and $(A_u(VQV^*),\s{\del}{VQV^*})$
for any unitary $m\times m$ matrix $V$. Therefore we can assume that $Q$ is a
diagonal strictly positive matrix: $Q=\mathrm{diag}(q_1,\ldots,q_m)>0$.

Let us fix a strictly positive matrix $Q=\mathrm{diag}(q_1,\ldots,q_m)$ and
let $\pi$ be the quotient map from $(A_u(Q),\s{\del}{Q})$ to its canonical Kac
quotient. By examining the commutation relations
\refeq{unirel} in the same way as those for $S_qU(2)$ we arrive at the
conclusion that $\pi(u_{k,l})=0$ whenever $q_k\neq q_l$ ($k,l=1,\ldots,m$).
Let $\s{B}{Q}$ be the universal $\cst$-algebra generated by $m^2$ elements of
a matrix $u=(u_{k,l})$ subject to relations \refeq{unirel} supplemented by
\[
u_{k,l}=0\text{ if }q_k\neq q_l.
\]
Now it is easy to see that $\s{B}{Q}$ is the free product of algebras
$\bigl(A_u(I_{n_p})\bigr)_{p=1,\ldots,N}$, where $n_1,\ldots,n_N$ are
multiplicities of distinct eigenvalues of $Q$. Also $\s{B}{Q}$ has a compact
quantum group structure of this free product (cf.~\cite[Theorem 1.1]{free}).
Since the obvious map from $A_u(Q)$ to $\s{B}{Q}$ is a compact quantum group
morphism, we are in the situation of Lemma \ref{klasy} (with $B^1=A_u(Q)$ and
$B^2=\s{B}{Q}$).

In the special case when all eigenvalues of $Q$ are pairwise different, the
$\cst$-algebra $\s{B}{Q}$ is isomorphic to $\cst(\FF_m)$. This algebra, in
turn, is residually finite dimensional by a theorem of Choi (\cite{choi}, see
also \cite[Corollary 4.5]{bl}). Therefore by the reasoning presented before
Corollary \ref{RDFkac}, $\s{B}{Q}=\cst(\FF_m)$ with its usual compact quantum
group structure is the canonical Kac quotient of $(A_u(Q),\s{\del}{Q})$.

In general it is not known whether the algebras $A_u(I_{n_p})$ are residually
finite dimensional. If this were the case then their free product $\s{B}{Q}$
would also be residually finite dimensional (\cite[Theorem 3.2]{el})
and we would conclude that
$\s{B}{Q}$ with its compact quantum group structure is the canonical Kac
quotient of $(A_u(Q),\s{\del}{Q})$. Nevertheless, even without the knowledge
of the structure $\s{B}{Q}$, Lemma \ref{klasy} tells us that the canonical Kac
quotient of $(A_u(Q),\s{\del}{Q})$ depends up to isomorphism only on the
quantum group $A_u(I_{n_1})*\cdots*A_u(I_{n_N})$ (with appropriate
comultiplication).

Now Corollary \ref{isoKac} says that the isomorphism class of
$(\s{\AP}{Q},\s{\del}{Q})$ depends only on the multiplicities of distinct
eigenvalues of $Q$.

The infinite dimensionality of $\s{\AP}{Q}$ follows from the fact that
$\s{\AP}{Q}$ contains the algebra of continuous functions on the Bohr
compactification of the free group $\FF_m$, which is an infinite group, as it
contains the group $\FF_m$ itself (by residual finite dimensionality of
$\cst(\FF_m)$). More precisely: the free group $\cst$-algebra $\cst(\FF_m)$ is
a quotient of the $\cst$-algebra of functions on the canonical Kac quotient of
$(A_u(Q),\s{\del}{Q})$. Therefore the commutative algebra $\Ci{\FF_m}$ embedds
into $\M{\Auhat(Q)}$ with a quantum group morphism. This morphism extends to
$\M{\Ci{\FF_m}}=\Cb{\FF_m}$. It is easy to see now that with this embedding all
almost periodic functions on $\FF_m$ are contained in $\s{\AP}{Q}$. In other
words $\s{\AP}{Q}$ contains the $\cst$ algebra of continuous functions on the
Bohr compactification of $\FF_m$ (cf.~Proposition \ref{almper}).

In particular, since the Bohr compactification of the free group is
not metrizable, the algebra of functions on this compactification is not
separable. It follows that $\s{\AP}{Q}$ is not separable.

\paragraph{MAP-quantum groups.}
A discrete group $\Gamma$ is said to be {\em maximally almost periodic}\/ or
an {\em MAP-group}\/ if the set of almost periodic functions on $\Gamma$
separates points of $\Gamma$. This is equivalent to saying that the canonical
homomorphism of $\Gamma$ into its Bohr compactification is injective. It is easy
to give examples of groups which are not maximally almost periodic (see
e.g.~\cite{vw}).
It is known, however, that a group $\Gamma$ is maximally almost
periodic if the $\cst$-algebra $\cst(\Gamma)$ is residually finite dimensional
(\cite[Remark 4.2(iii)]{bl}). It turns out that we can generalize this result to
$\cst$-algebraic discrete quantum groups.

In what follows we shall call a
$\cst$-algebraic discrete quantum group $(A,\s{\del}{A})$ {\em maximally
almost periodic}\/ if the range of the canonical map $\chi$ from its Bohr
compactification to $(A,\s{\del}{A})$ is strictly dense in $\M{A}$. This
condition clearly corresponds to the injectivity of the homomophism from a
discrete group into its Bohr compactification.
\pagebreak

\begin{proposition}\label{mapy}
\noindent
\begin{itemize}
\item[{\rm (1)}]
Let $(A,\s{\del}{A})$ be a $\cst$-algebraic discrete quantum group such that its
universal dual $(B,\s{\del}{B})$ has the property that $B$ is residually finite
dimensional. Then $(A,\s{\del}{A})$ is unimodular and maximally almost
periodic.
\item[{\rm (2)}] Any maximally almost periodic $\cst$-algebraic discrete
quantum group is unimodular.
\end{itemize}
\end{proposition}

\proof
{\sc Ad (1).}
Let $(\AP,\s{\del}{A})$ be the Bohr compactification of $(A,\s{\del}{A})$.
The unimodularity of $(A,\s{\del}{A})$ follows from Corollary \ref{RDFkac} and
Theorem \ref{unim}. The algebra $\AP$ is generated by matrix elements of finite
dimensional unitary representations of $(A,\s{\del}{A})$. Every such
representation is of the form $(\id\tens\pi)u^*$, where $u$ is the universal
bicharacter for the duality of $(A,\s{\del}{A})$ and $(B,\s{\del}{B})$
and $\pi$ is a finite dimensional representation of $B$ (cf.~Section
\ref{dowody}). Let $\ph$ be a non zero continuous functional on $A$. Then there
exists an $x\in\AP$ such that $\ph(x)$ is non zero (we use here the canonical
extension of continuous functionals on $A$ to $\M{A}$). Indeed: if $\ph(x)=0$
for all $x$ of the form $x=(\id\tens\pi)u^*$ with $\pi$ a finite dimensional
representation of $B$ then the element $(\ph\tens\id)u^*$ is in the intersection
of kernels of all finite dimensional representations of $B$. It follows that
$(\ph\tens\id)u^*=0$ and consequently $\ph=0$. Since continuous functionals on
$A$ are precisely the strictly continuous functionals on $\M{A}$ and $\M{A}$ is
the strict completion of $A$, it follows that $\AP$ is strictly dense in
$\M{A}$.

{\sc Ad(2).}
Now assume that $(A,\s{\del}{A})$ is maximally almost periodic. Then the
coinverse of the compact quantum group of Kac type $(\AP,\s{\del}{A})$ is
bounded and $\AP$ is strictly dense in $\M{A}$. The canonical extension of
$\s{\kap}{A}$ to $\M{A}$ coincides with its strict closure. The closure of a map
which is bounded on a dense set is bounded. In particular the restriction of
$\s{\kap}{A}$ to $A$ is bounded and $(A,\s{\del}{A})$ is unimodular by Theorem
\ref{unim}.
\qed

The second part of Proposition \ref{mapy} describes the extent to which Bohr
compactification remembers the original quantum group. The map from a group to
its Bohr compactification cannot be injective for non unimodular discrete
quantum groups.

\section*{Acknowledgments}
The author is deeply indebted to Stefaan Vaes for helpful suggestions and
many fruitful discussions on the subject of Bohr compactification.
The paper was written during author's stay at Mathematisches Institut der
Wesf\"alischen Wilhelms-Universit\"at in M\"unster. He wishes to thank
professor Joachim Cuntz and colleagues from the institute for creating
unrivaled conditions for scientific activity.


\begin{thebibliography}{XX}
\bibitem{banica}
{\sc Banica, T.:} Le groupe quantique compact libre $\mathrm{U}(n)$, {\it
Comm.~Math.~Phys.} \textbf{190} no.~1 (1997), 143--172.
\bibitem{bl}
{\sc Bekka, M.B.~\& Louvet, N.:} Some properties of $\mathrm{C}^*$-algebras
associated to discrete linear groups. In: {\it $\mathrm{C}^*$-Algebras, Proc.~of
the SFB-Workshop, J.~Cuntz, S.~Etcherhoff (Eds.), M\"unster 1999,}\/ pp.~1--22,
Springer 1999.
\bibitem{choi}
{\sc Choi, M.D.:} The full $\mathrm{C}^*$-algebra of the free group on two
generators. {\it Pac.~J.~Math.}\/ \textbf{87} no.~1 (1980), 41--48.
\bibitem{ev}
{\sc Enock, M.~\& Vallin, J.-M.:} $\mathrm{C}^*$-alg\`ebres de Kac et
alg\`ebres de Kac. {\it Proc.~London Math.~Soc.} \textbf{66} no.~3 (1993),
619--650.
\bibitem{el}
{\sc Exel, R.~\& Lorring, T.A.:} Finite dimensional representations of free
product $\mathrm{C}^*$-algebras. {\it Int.~J.~Math.}\/ \textbf{3} no.~4
(1992), 469--476.
\bibitem{hm}
{\sc Hajac, P.M.~\& Masuda, T.:} Quantum double torus. {\it
C.~R.~Acad.~Sci.~Paris Ser.~I Math.} \textbf{327} no.~6 (1998), 553--558.
\bibitem{holm}
{\sc Holm, P.:} On the Bohr compactification. {\it Math.~Ann.}\/
\textbf{156} (1964), 34--46.
\bibitem{loom}
{\sc Loomis, L.H.:} {\it An introduction to abstract harmonic
analysis.}\/ Van Nostrad, Princeton N.J.~1953.
\bibitem{mvd}
{\sc Maes, A.~\& Van Daele, A.:} Notes on compact quantum groups. {\it
Nieuw Arch.~Wisk.~(4)}\/ \textbf{16} nos.~1--2 (1998), 73--112.
\bibitem{pw}
{\sc Podle\'{s}, P.~\& Woronowicz, S.L.:} Quantum deformation of
Lorentz group. {\it Comm. Math.~Phys.}\/ \textbf{130} no.~2 (1990), 381--431.
\bibitem{cdqg}
{\sc So\l{}tan, P.M.:} Compactifications of discrete quantum groups.
Preprint des SFB 478 des Mathematisches Institut der
Wesf\"alischen Wilhelms-Universit\"at M\"unster (2003). To appear
in {\it Algebr.~Represent.~Theory.}
\bibitem{mha}
{\sc Van Daele, A.:} Multiplier Hopf algebras. {\it Trans.~AMS}\/
\textbf{342} no.~2 (1994), 917--932.
\bibitem{dqg}
{\sc Van Daele, A.:} Discrete quantum groups. {\it J.~Algebra}\/
\textbf{180} no.~2 (1996), 431--444.
\bibitem{afgd}
{\sc Van Daele, A.:} An algebraic framework for group duality. {\it
Adv.~in Math.}\/ \textbf{140} no.~2 (1998), 323--366.
\bibitem{vawa}
{\sc Van Daele, A.~\& Wang, S.:} Universal quantum groups. {\it Int.~J.~Math.}
\textbf{7} no.~2 (1996), 747--764.
\bibitem{va}
{\sc Vaes, S.} Personal communication.
\bibitem{vw}
{\sc Von Neumann, J.~\& Wigner, E.~P.:} Minimally almost periodic groups.
{\it Ann. of Math.~(2)}\/ \textbf{41} (1940), 746--750.
\bibitem{free}
{\sc Wang, S.:} Free products of compact quantum groups. {\it
Comm.~Math.~Phys.} \textbf{167} no.~3 (1995), 671--692.
\bibitem{iso}
{\sc Wang, S.:} Structure and isomorphism classes of compact quantum groups
$A_u(Q)$ and $B_u(Q)$. {\it J.~Op.~Th.} \textbf{167} no.~3, suppl.~(2002),
573--583.
\bibitem{weil}
{\sc Weil, A.:} {\it L'integration dans les groupes topologiques et
ses applications.}\/ Herman, Paris 1965.
\bibitem{pseudo}
{\sc Woronowicz, S.L.:} Pseudospaces, pseudogroups and Pontryagin
duality. In: {\it Mathematical problems in theoretical Physics,
Proc.~Int.~Conf.~Math.~Phys., Lausanne 1979,}\/ pp.~407--422, Lecture
Notes in Phys.~\textbf{116}, Springer Berlin--Heidelberg--New York
1980.
\bibitem{SUq2}
{\sc Woronowicz, S.L.:} Twisted $SU(2)$ groups. An example of
non-commutative differential calculus. {\it Publ.~RIMS Kyoto Univ.}
\textbf{23} no.~1 (1987), 117-181.
\bibitem{pseudogr}
{\sc Woronowicz, S.L.:} Compact matrix pseudogroups. {\it
Comm.~Math.~Phys.}\/\textbf{111} no.~4 (1987), 613--665.
\bibitem{unbo}
{\sc Woronowicz, S.L.:} Unbounded elements affiliated with
$\mathrm{C}^*$-algebras and non-compact quantum groups. {\it
Comm.~Math.~Phys.}\/ 136 no.~2 (1991), 399--432.
\bibitem{gen}
{\sc Woronowicz, S.L.:} $\mathrm{C}^*$-algebras generated by unbounded
elements. {\it Rev. Math.~Phys.}\/ \textbf{7} no.~3 (1995), 481--521.
\bibitem{cqg}
\sloppy
{\sc Woronowicz, S.L.:} Compact quantum groups. In: {\it Sym\'etries
quantiques, les Houches, Session LXIV 1995,}\/ pp. 845--884, Elsevier 1998.
\bibitem{wn}
{\sc Woronowicz, S.L.~\& Napi\'orkowski, K.:} Operator theory in
the $\mathrm{C}^*$-algebra framework. {\it Rep.~Math.~Phys.}\/
\textbf{31} no.~3 (1992), 353--371.
\end{thebibliography}
\end{document}